\documentclass[reqno]{amsart}
\usepackage{hyperref}
\usepackage{amsmath}
\begin{document}

\title[\hfilneg \hfil A note on the Uniqueness of entire functions  sharing two pairs of values with its  difference operator]
{A note on the Uniqueness of entire functions  sharing two pairs of values with its  difference operator}

\author[XiaoHuang Huang \hfil \hfilneg]
{XiaoHuang Huang}

\address{XiaoHuang Huang: Corresponding author\newline
Department of Mathematics,  Department of Mathematics, Sun Yat-sen University, Guangzhou, 510275, P.R. China}
\email{1838394005@qq.com}

\begin{abstract}
In this paper, we investigate the sharing values problem that entire function $f(z)$ and its first order difference operator $\Delta_{\eta}f(z)$ share two distinct pairs of finite values IM.  We  prove: Let $f(z)$ be a non-constant entire function of hyper-order less than $1$, let $\eta$ be a non-zero complex number, and let $a$ be a nonzero finite number. Then there exists no such entire function so that  $ f(z)$ and $\Delta_{\eta}f(z)$ share $(0,0)$  and $(a,-a)$ IM. Furthermore, using a result in Wang-Chen-Hu \cite{wch}, we obtain some uniqueness results that when $ f(z)$ and $\Delta_{\eta}f(z)$ share $a\neq0$  and $-a$ IM.
\end{abstract}

\maketitle
\numberwithin{equation}{section}
\newtheorem{theorem}{Theorem}[section]
\newtheorem{lemma}[theorem]{Lemma}
\newtheorem{remark}[theorem]{Remark}
\newtheorem{corollary}[theorem]{Corollary}
\newtheorem{example}[theorem]{Example}
\newtheorem{problem}[theorem]{Problem}
\allowdisplaybreaks

\section{ main results}
Let $\Bbb C$ denote the complex plane and $f(z)$ a meromorphic function on $\Bbb C$. In this paper, we assume that the reader is familiar with the fundamental results
and the standard notation of the Nevanlinna value distribution theory, see(\cite{h3,y1,y2}). In addition, $S(r, f) = o(T(r, f))$, as $r\to \infty $ outside of a possible exceptional set of finite logarithmic measure. Define
 $$\rho(f)=\varlimsup_{r\rightarrow\infty}\frac{log^{+}T(r,f)}{logr},$$
 as the  order  of $f$.

Let $a,b$ be  two complex numbers. We say that two non-constant meromorphic functions $f(z)$ and $g(z)$ share a pair  $(a,b)$ IM (CM) if $f(z)-a$ and $g(z)-b$ have the same zeros ignoring multiplicities (counting multiplicities). Moreover, we introduce the following denotations : $\overline{N}_{(m,n)}(r,\frac{1}{f-a})$ denotes the reduced counting function of both zeros $f(z)-a$ with multiplicity $m$ and zeros of $g(z)-b$ with multiplicity $n$. $\overline{N}_{[n]}(r,\frac{1}{f(z)-a})$ denotes the reduced counting function of all zeros of $f(z)-a$ with multiplicity $n$.

Let $a$ be  a complex numbers. We say that two non-constant meromorphic functions $f(z)$ and $g(z)$  share value $a$ CM almost if
$$N(r,\frac{1}{f-a})+2N(r,\frac{1}{g-a})-N(r,f=a=g)=S(r,f),$$
where $N(r,f=a=g)$ means the common zeros of $f-a$ and $g-a$.

Let $f(z)$ be a non-constant meromorphic function, and let $\eta$ be a nonzero complex constant. We define the difference operator
$$\Delta_{\eta}f(z)=f(z+\eta)-f(z).$$

In 1977, Rubel and Yang \cite{ruy}  considered the uniqueness of an entire function and its derivative. They proved.

\

{\bf Theorem 1} \ Let $f(z)$ be a transcendental entire function, and let $a, b$ be two finite distinct complex values. If $f(z)$ and $f'(z)$
 share $a, b$ CM, then $f(z)\equiv f'(z)$.

Mues and Steinmetz \cite {muess} improved Theorem 1 and proved

\

{\bf Theorem 2} \ Let $f(z)$ be a transcendental entire function, and let $a, b$ be two finite distinct complex values. If $f(z)$ and $f'(z)$
 share $a, b$ IM, then $f(z)\equiv f'(z)$.

Recently, the difference analogue of the lemma on the logarithmic derivative and Nevanlinna theory for the difference operator have been founded, which bring about a number of papers $[1-6, 8-14, 19]$ focusing on the uniqueness study of meromorphic functions sharing some values with their difference operators.  It is well known that $\Delta_{\eta}f(z)$ can be considered as the difference counterpart of $f'(z)$ in Theorem 1.

Corresponding to Theorem 1, Chen-Yi \cite{cy} proved

 \

{\bf Theorem 3}
 Let $f(z)$ be a transcendental entire function of finite non-integer order,   let $\eta$ be a non-zero complex number  and let $a$ and $b$ be two distinct complex values. If $ f(z)$ and $\Delta_{\eta}f(z)$ share $a$, $b$ CM, then $ f(z)\equiv \Delta_{\eta}f(z)$.

In \cite{cy}, the authors conjectured that it's not necessary that the order is not a positive integer. In 2014, Zhang-Liao \cite{zl},  Liu-Yang-Fang \cite{lyf} confirmed  the conjecture. They proved

\

{\bf Theorem 4}
 Let $f(z)$ be a transcendental entire function of finite order,  let $\eta$ be a non-zero complex number, and let $a, b$ be two finite distinct complex values. If $ f(z)$ and $\Delta_{\eta}f(z)$ share $a$, $b$ CM, then $ f(z)\equiv \Delta_{\eta}f(z)$.

Corresponding to Theorem 3, Cui-Chen \cite{cc}, Li-Yi-Kang \cite{lyk}, L\"u-L\"u \cite{ll} proved

\

{\bf Theorem 5}
 Let $f(z)$ be a transcendental meromorphic function of finite order, let $\eta$ be a non-zero complex number, and let $a$ and $b$ be two distinct complex numbers. If $f(z)$ and $\Delta_{\eta}f(z)$ share $a$, $b$, $\infty$ CM, then $f(z)\equiv \Delta_{\eta}f(z)$.

By Theorem 1-Theorem 5, it's natural to pose the following question.\\

\

{\bf Question 1} Let $f(z)$ be a transcendental entire function of finite order, let $\eta$ be a non-zero complex number, and let $a$ and $b$ be two distinct complex numbers. If $ f(z)$ and $\Delta_{\eta}f(z)$ share $a$ and $b$ IM, is $ f(z)\equiv \Delta_{\eta}f(z)$?\\

In 2017, Li-Duan-Chen \cite{ldc} proved

\

{\bf Theorem 6}
 Let $f(z)$ be a  transcendental entire function of finite order, let $\eta$ be a non-zero complex number,  and let $a$ be a nonzero complex number. If $f(z)$ and $\Delta_{\eta}f(z)$ share $0$ CM and share $a$ IM, then $f(z)\equiv \Delta_{\eta}f(z)$.

In \cite{ldc}, the authors posed the following question.

\

{\bf Question 2} Let $f(z)$ be a transcendental entire function of finite order, let $\eta$ be a non-zero complex number, and let $a$ and $b$ be two distinct complex numbers. If $ f(z)$ and $\Delta_{\eta}f(z)$ share $a$ CM and share $b$ IM, is $ f(z)\equiv \Delta_{\eta}f(z)$?

In 2021\cite{h4}, we give a positive answer to  Question 2. We prove.

\

{\bf Theorem 7} Let $f(z)$ be a transcendental entire function of finite order, let $\eta$ be a non-zero complex number, and let $a$ and $b$ be two distinct complex numbers. If $ f(z)$ and $\Delta_{\eta}f(z)$ share $a$ CM and share $b$ IM, then $ f(z)\equiv \Delta_{\eta}f(z)$.

In this paper, we study uniqueness of entire function sharing two pairs of finite values with its first difference operator. We obtain the following result.

\

{\bf Theorem 8}  Let $f(z)$ be a nonconstant entire function of finite order,  $\eta$  a nonzero complex number, and $i=1,2$. Let $a_{i}$ and $b_{i}$ be four finite complex numbers such that $a_{i}\neq b_{i}$, and let
$$\varphi(z)=\frac{f'(z)((a_{2}-b_{2})f(z)-(a_{1}-b_{1})\Delta_{\eta}f(z)- a_{2}b_{1}+a_{1}b_{2})}{(f(z)-a_{1})(f(z)-b_{1})},$$
If $f(z)$ and $\Delta_{\eta}f(z)$ share $(a_{1},a_{2})$, $(b_{1},b_{2})$ IM,  then either 
$$(a_{2}-b_{2})f(z)-(a_{1}-b_{1})\Delta_{\eta}f(z)\equiv a_{2}b_{1}-a_{1}b_{2},$$
or $\varphi(z)$ is a periodic entire function with period $\eta$.

By Theorem 8, we can partially solve Question 1 in the case of $(0,0)$ and $(a,-a)$. We prove.

\

{\bf Theorem 9} Let $f(z)$ be a non-constant entire function of finite order, let $\eta$ be a non-zero complex number, and let $a$ be a nonzero finite number. Then there exists no such entire function so that  $ f(z)$ and $\Delta_{\eta}f(z)$ share $(0,0)$  and $(a,-a)$ IM.

It is natural to ask whether $(0,0)$ can be replaced by $(b,-b)$, where $a\neq b$.

\

{\bf Question 2} Let $f(z)$ be a non-constant entire function of finite order, let $\eta$ be a non-zero complex number, and let $a,b$ be two distinct finite number. Then there exists no such entire function so that  $ f(z)$ and $\Delta_{\eta}f(z)$ share $(a,-a)$  and $(b,-b)$ IM.

Furthermore, as an application of Theorem 8, we obtain a uniqueness result using Lemma 2.9 in  Section 2. We obtain.

\

{\bf Theorem 10} Let $f(z)$ be a non-constant entire function of $\rho(f)<1$, let $\eta$ be a non-zero complex number, and let $a$ be a nonzero finite number. if  $ f(z)$ and $\Delta_{\eta}f(z)$ share $a$  and $-a$ IM, then $f(z)\equiv\Delta_{\eta}f(z)$.

\section{Some Lemmas}
\begin{lemma}\label{211} Let $f$ be a non-constant meromorphic function with finite order,  and let $c$ be a non-zero complex number. Then
$$m(r,\frac{f(z+c)}{f(z)})=S(r, f),$$
for all r outside of a possible exceptional set E with finite logarithmic measure.
\end{lemma}

\begin{lemma}\label{22l}\cite{h3,y1,y2}
Suppose $f_{1}(z),f_{2}(z)$ are two non-constant meromorphic functions in the complex plane, then
$$N(r,f_{1}f_{2})-N(r,\frac{1}{f_{1}f_{2}})=N(r,f_{1})+N(r,f_{2})-N(r,\frac{1}{f_{1}})-N(r,\frac{1}{f_{2}}).$$
\end{lemma}

\begin{lemma}\label{23l}\cite {h2}
Let $f(z)$ be a non-constant meromorphic function of finite order, and let $\eta\neq0$ be a finite complex number. Then
$$T(r,f(z+\eta))=T(r, f(z))+S(r,f).$$
\end{lemma}

\begin{lemma}\label{24l}\cite{h3,y2} Let $f(z)$ be a nonconstant meromorphic function, and let $P(f)=a_{0}f^{p}(z)+a_{1}f^{p-1}(z)+\cdots+a_{p}(a_{0}\neq0)$ is a polynomial of degree $p$ with constant coefficients $a_{j}(j=0,1,\ldots,p)$.Suppose that $b_{j}(j=0,1,\ldots,q)(q>p)$. Then
$$m(r,\frac{P(f)f'(z)}{(f(z)-b_{1})(f(z)-b_{2})\cdots(f(z)-b_{q})})=S(r,f).$$
\end{lemma}

\begin{lemma}\label{251}\cite{h3,y2} Suppose that $f(z)$ is a meromorphic function in the complex plane and $p(f)= a_{0}f^{n}(z)+a_{1}f^{n-1}(z)+\cdots +a_{n}$ , where $a_{0}(\not\equiv0)$, $a_{1}$,$\cdots$,$a_{n}$ are complex numbers. Then
$$T(r,p(f))= nT(r,f(z))+S(r,f).$$
\end{lemma}
In 2020, Wang-Chen-Hu \cite{wch} consider the  entire  solutions of the following quadratic functional equation:
\begin{align}
&\{a_{0}f(z)+a_{1}f'(z)+a_{2}f(z+\eta)\}^{2}+2a_{3}\{a_{0}f(z)+a_{1}f'(z)+a_{2}f(z+\eta)\}\{b_{0}f(z)\notag\\
&+b_{1}f'(z)+b_{2}f(z+\eta)\}+\{b_{0}f(z)+b_{1}f'(z)+b_{2}f(z+\eta)\}^{2}=e^{c_{1}z+c_{2}},
\end{align}
where $f(z)$ is an entire function of finite order, and $a_{i},b_{i}(i=0,1,2)$ and $c_{j}(j=1,2)$ are finite complex number.

Let  six complex numbers $a_{i},b_{i}(i=0,1,2)$ be a matrix $A$ such that 
\begin{equation} 
rank\left ( 
\begin{array}{ccc} 
a_{0} & a_{1} & a_{2}\\ 
b_{0} & b_{1} & b_{2}\\ 
\end{array}
\right)=2.
\end{equation}
They obtain the following result.
\begin{lemma}\label{261}\cite{wch} Take complex number $a_{i},b_{i}(i=0,1,2)$, $c_{j}(j=1,2)$, $\eta$, $a_{3}$ with $\eta\neq0$, $a_{i}^{2}\neq1,0$ and ((2.2) holds. If the equation (2.1) 
has entire solutions $ f$ with finite order, then it only has the following form:
\begin{align}
f(z)=C_{1}e^{dz}+C_{2}e^{-dz}+C_{3}e^{C_{0}z},
\end{align}
where $d\neq0$, $C_{k}(k=0,1,2,3)$ are constants. Moreover $C_{1}C_{2}\neq0$, $C_{3}=0$ and $C_{0}=0$ if $a_{0}b_{1}-a_{1}b_{0}\neq0$, and $a_{1}b_{2}-a_{2}b_{1}\neq0$.\end{lemma}

\begin{lemma}\label{271}\cite{h3,y1}
Let $f$ and $g$ be two non-constant meromorphic functions with $\lambda(f)$ and $\lambda(g)$ as their orders respectively. Then $$\lambda(fg)\leq max\{\lambda(f),\lambda(g)\}.$$
\end{lemma}

\begin{lemma}\label{281}\cite{fw}
Let $f$  be a meromorphic functions, let $\eta\neq0$ be a finite complex number, and let $n$ be a positive integer. If $)\Delta_{\eta}^{n}f\equiv0$, then either $\rho(f)\geq1$ or $f$ is a polynomial with $\deg(f)\leq n-1$.
\end{lemma}

\section{The proof of Theorem 8 }
We prove Theorem 8 by contradiction. Suppose that $(a_{2}-b_{2})f-(a_{1}-b_{1})\Delta_{\eta}f- a_{2}b_{1}+a_{1}b_{2}\not\equiv0$.  Since $f$ is an  entire function of finite order, and that $f$ and $\Delta_{\eta}f$ share $(a_{1},a_{2})$ and $(b_{1},b_{2})$ IM, then  by the First and  Second Fundamental Theorem and Lemma 2.1, we get
\begin{eqnarray*}
\begin{aligned}
T(r,f)&\leq \overline{N}(r,\frac{1}{f-a_{1}})+\overline{N}(r,\frac{1}{f-b_{1}})+S(r,f)\\
&= \overline{N}(r,\frac{1}{\Delta_{\eta}f-a_{2}})+\overline{N}(r,\frac{1}{\Delta_{\eta}f-b_{2}})+S(r,f)\\
&\leq N(r,\frac{1}{(a_{2}-b_{2})f-(a_{1}-b_{1})\Delta_{\eta}f- a_{2}b_{1}+a_{1}b_{2}})+S(r,f)\\
&\leq T(r,(a_{2}-b_{2})f-(a_{1}-b_{1})\Delta_{\eta}f- a_{2}b_{1}+a_{1}b_{2})+S(r,f)\\
&= m(r,(a_{2}-b_{2})f-(a_{1}-b_{1})\Delta_{\eta}f)+S(r,f)\\
&\leq m(r,f)+m(r,1-\frac{\Delta_{\eta}f}{f})+S(r,f)\\
&\leq T(r,f)+S(r,f).
\end{aligned}
\end{eqnarray*}

Thus we have
\begin{eqnarray}
T(r,f)=\overline{N}(r,\frac{1}{f-a_{1}})+\overline{N}(r,\frac{1}{f-b_{1}})+S(r,f).
\end{eqnarray}

Set
\begin{eqnarray}
\varphi=\frac{f'((a_{2}-b_{2})f-(a_{1}-b_{1})\Delta_{\eta}f- a_{2}b_{1}+a_{1}b_{2})}{(f-a_{1})(f-b_{1})}\\
\psi=\frac{\Delta_{\eta}f'((a_{2}-b_{2})f-(a_{1}-b_{1})\Delta_{\eta}f- a_{2}b_{1}+a_{1}b_{2})}{(\Delta_{\eta}f-a_{2})(\Delta_{\eta}f(z)-b_{2})}.
\end{eqnarray}

Noting that $f$ is an entire function of finite order and that $f$ and $\Delta_{\eta}f$ share $(a_{1},a_{2})$ and $(b_{1},b_{2})$ IM,  we know that $\varphi$ is an entire function by (3.2). By Lemma 2.1 and Lemma 2.3,  we have
\begin{eqnarray*}
\begin{aligned}
&T(r,\varphi)=m(r,\varphi(z))=m(r,\frac{f'((a_{2}-b_{2})f-(a_{1}-b_{1})\Delta_{\eta}f- a_{2}b_{1}+a_{1}b_{2})}{(f-a_{1})(f-b_{1})})+S(r,f)\notag\\
&\leq m(r,\frac{f'f}{(f-a_{1})(f-b_{1})})+m(r,1-\frac{f'\Delta_{\eta}f}{(f-a_{1})(f-b_{1})})\notag\\
&m(r,\frac{f'(a_{2}b_{1}-a_{1}b_{2})}{(f-a_{1})(f-b_{1})})+S(r,f)=S(r,f),
\end{aligned}
\end{eqnarray*}
that is
\begin{align}
T(r,\varphi)=S(r,f).
\end{align}
Let $d_{i}= a_{i}+k(a_{i}-b_{i})$ be a finite value, where $k\neq0$ and $i=1,2$, then it follows from Nevanlinna’s second fundamental theorem and (3.1) that
\begin{align}
 2T(r,f)&\leq \overline{N}(r,\frac{1}{f-a_{1}})+\overline{N}(r,\frac{1}{f-b_{1}})+\overline{N}(r,\frac{1}{f-d_{1}})+S(r,f)\notag\\
 &=T(r,f)+\overline{N}(r,\frac{1}{f-d_{1}})+S(r,f),
  \end{align}
which implies that
\begin{align}
m(r,\frac{1}{f-d_{1}})=S(r,f).
\end{align}By (3.2) we get
\begin{align}
\varphi f^{2}=f'((a_{2}-b_{2})f-(a_{1}-b_{1})\Delta_{\eta}f)+(a_{1}+b_{1})\varphi f+a_{1}b_{1}\varphi+(a_{1}b_{2}-a_{2}b_{1})f',
\end{align}
and
\begin{align}
\varphi_{\eta} f^{2}_{\eta}=f'_{\eta}((a_{2}-b_{2})f_{\eta}-(a_{1}-b_{1})\Delta_{\eta}f_{\eta})+(a_{1}+b_{1})\varphi_{\eta} f_{\eta}+a_{1}b_{1}\varphi_{\eta}+(a_{1}b_{2}-a_{2}b_{1})f'_{\eta}.
\end{align}
By (3.7) and (3.8), we get
\begin{align}
&\Delta_{\eta}\varphi f^{2}=f'_{\eta}((a_{2}-b_{2})\Delta_{\eta}f-(a_{1}-b_{1})\Delta_{\eta}^{2}f)+\Delta_{\eta}f'((a_{2}-b_{2})f-(a_{1}-b_{1})\Delta_{\eta}f)\notag\\
&+(a_{1}+b_{1})[\varphi_{\eta}f_{\eta}-\varphi f]-\varphi_{\eta}[f^{2}_{\eta}-f^{2}]-a_{1}b_{1}\Delta_{\eta}\varphi+(a_{1}b_{2}-a_{2}b_{1})\Delta_{\eta}f'.
\end{align}

{\bf Case 1.}\quad $\Delta_{\eta}\varphi\not\equiv0$.

By (3.9), Lemma 2.1 and Lemma 2.2, we  obtain
\begin{eqnarray*}
\begin{aligned}
&2m(r,f)=m(r,f'_{\eta}((a_{2}-b_{2})\Delta_{\eta}f-(a_{1}-b_{1})\Delta_{\eta}^{2}f)+\Delta_{\eta}f'((a_{2}-b_{2})f-(a_{1}-b_{1})\Delta_{\eta}f)\notag\\
&+(a_{1}+b_{1})[\varphi_{\eta}f_{\eta}-\varphi f]-\varphi_{\eta}[f^{2}_{\eta}-f^{2}]+(a_{1}b_{2}-a_{2}b_{1})\Delta_{\eta}f')+S(r,f)\notag\\
&\leq m(r,f)+m(r,\frac{f'_{\eta}((a_{2}-b_{2})\Delta_{\eta}f-(a_{1}-b_{1})\Delta_{\eta}^{2}f)+\Delta_{\eta}f'((a_{2}-b_{2})f-(a_{1}-b_{1})\Delta_{\eta}f)}{f}\\
&-\frac{\varphi_{\eta}[f^{2}_{\eta}-f^{2}]}{f})+m(r,\frac{(a_{1}+b_{1})[\varphi_{\eta}f_{\eta}-\varphi f]+(a_{1}b_{2}-a_{2}b_{1})\Delta_{\eta}f'}{f})+S(r,f)\\
&\leq m(r,\frac{f'_{\eta}((a_{2}-b_{2})\Delta_{\eta}f-(a_{1}-b_{1})\Delta_{\eta}^{2}f)+\Delta_{\eta}f'((a_{2}-b_{2})f-(a_{1}-b_{1})\Delta_{\eta}f)}{f\Delta_{\eta}f}\\
&-\frac{\varphi_{\eta}[f^{2}_{\eta}-f^{2}]}{f\Delta_{\eta}f})+m(r,f)+m(r,\Delta_{\eta}f)+S(r,f)\\
&\leq m(r,f)+m(r,\Delta_{\eta}f)+m(r,\frac{f'_{\eta}((a_{2}-b_{2})\Delta_{\eta}f-(a_{1}-b_{1})\Delta_{\eta}^{2}f)}{f\Delta_{\eta}f})\\
&+m(r,\frac{\Delta_{\eta}f'((a_{2}-b_{2})f-(a_{1}-b_{1})\Delta_{\eta}f)}{f\Delta_{\eta}f})+m(r,\frac{\varphi_{\eta}[f^{2}_{\eta}-f^{2}]}{f\Delta_{\eta}f})+S(r,f).\\
&\leq m(r,f)+m(r,\Delta_{\eta}f)+S(r,f).
\end{aligned}
\end{eqnarray*}

It follows
\begin{align}
T(r,f)\leq T(r,\Delta_{\eta}f)+S(r,f).
\end{align}

Obviously,
 \begin{align}
T(r,\Delta_{\eta}f)\leq T(r,f)+S(r,f).
\end{align}
It follows from (3.10) and (3.11) that
\begin{align}
T(r,f)= T(r,\Delta_{\eta}f)+S(r,f).
\end{align}
By the First and Second  Fundamental Theorem, (3.1) and (3.12), we have
\begin{eqnarray*}
\begin{aligned}
2T(r,f)&\leq 2T(r,\Delta_{\eta}f)+S(r,f)\\
&\leq\overline{N}(r,\frac{1}{\Delta_{\eta}f-a_{2}})+\overline{N}(r,\frac{1}{\Delta_{\eta}f-b_{2}})+\overline{N}(r,\frac{1}{\Delta_{\eta}f-d_{2}})+S(r,f)\\
&\leq \overline{N}(r,\frac{1}{f-a_{1}})+\overline{N}(r,\frac{1}{f-b_{1}})+T(r,\frac{1}{\Delta_{\eta}f-d})-m(r,\frac{1}{\Delta_{\eta}f-d_{2}})\\
&+S(r,f)\leq T(r,f)+T(r,\Delta_{\eta}f)-m(r,\frac{1}{\Delta_{\eta}f-d_{2}})+S(r,f)\\
&\leq 2T(r,f)-m(r,\frac{1}{\Delta_{\eta}f-d_{2}})+S(r,f).
\end{aligned}
\end{eqnarray*}
Thus
\begin{eqnarray}
m(r,\frac{1}{\Delta_{\eta}f-d_{2}})=S(r,f).
\end{eqnarray}

By the First Fundamental Theorem, Lemma 2.1, Lemma 2.3, Lemma 2.4, (3.5), (3.6), (3.12), (3.13), and  $f(z)$ is an entire function of finite order, we obtain
\begin{eqnarray*}
\begin{aligned}
m(r,\frac{f-d_{1}}{\Delta_{\eta}f-d_{2}})&\leq T(r,\frac{f-d_{1}}{\Delta_{\eta}f-d_{2}})-N(r,\frac{f-d_{1}}{\Delta_{\eta}f-d_{2}})+S(r,f)\\
&=m(r,\frac{\Delta_{\eta}f-d_{2}}{f-d_{1}})+N(r,\frac{\Delta_{\eta}f-d_{2}}{f-d_{1}})-N(r,\frac{f-d_{1}}{\Delta_{\eta}f-d_{2}})+S(r,f)\\
&\leq N(r,\frac{1}{f-d_{1}})-N(r,\frac{1}{\Delta_{\eta}f-d_{2}})+S(r,f)\\
&=T(r,\frac{1}{f-d_{1}})-T(r,\frac{1}{\Delta_{\eta}f-d_{2}})+S(r,f)\\
&=T(r,f)-T(r,\Delta_{\eta}f)+S(r,f)=S(r,f),
\end{aligned}
\end{eqnarray*}
which implies that
\begin{eqnarray}
m(r,\frac{f-d_{1}}{\Delta_{\eta}f-d_{2}})=S(r,f).
\end{eqnarray}

By (3.3) we have
\begin{eqnarray}
\psi=[\frac{a_{2}-d_{2}}{a_{2}-b_{2}}\frac{\Delta_{\eta}f'}{\Delta_{\eta}f-a_{2}}-\frac{b_{2}-d_{2}}{a_{2}-b_{2}}\frac{\Delta_{\eta}f'}{\Delta_{\eta}f-b_{2}}][\frac{f-d_{1}}{\Delta_{\eta}f-d_{2}}-1].
\end{eqnarray}
Since $f$ is an entire function, and that $f$ and $\Delta_{\eta}f$ share $(a_{1},a_{2})$ and $(b_{1},b_{2})$ IM, we know that  $\psi$ is an entire function. Then by (3.14), (3.15), Lemma 2.1 and Lemma 2.2, we  get
\begin{align}
&T(r,\psi)=m(r,\psi)\leq m(r,\frac{a_{2}-d_{2}}{a_{2}-b_{2}}\frac{\Delta_{\eta}f'}{\Delta_{\eta}f-a_{2}})\notag\\
&+m(r,\frac{b_{2}-d_{2}}{a_{2}-b_{2}}\frac{\Delta_{\eta}f'}{\Delta_{\eta}f-b_{2}})+m(r,\frac{f-d_{1}}{\Delta_{\eta}f-d_{2}}-1)=S(r,f).
\end{align}

Now let $z_{1}$ be a zero of $f-a_{1}$ and $\Delta_{\eta}f-a_{2}$ with multiplicities $m$ and $n$, respectively. Using Taylor series expansions,  and by calculating we get $n\varphi(z_{1})-m\psi(z_{1})=0$. Let$$H_{n,m}=n\varphi-m\psi,$$
where $m$ and $n$ are positive integers.
Next, we consider two subcases.

{\bf Subcase 1.1} \quad $H_{n,m}\equiv0$ for some positive integers $m$ and $n$. That is $n\varphi\equiv m\psi$. Then  we have
$$n(\frac{f'(z)}{f(z)-a}-\frac{f'(z)}{f(z)-b})\equiv m(\frac{\Delta_{\eta}f'(z)}{\Delta_{\eta}f(z)-a}-\frac{\Delta_{\eta}f'(z)}{\Delta_{\eta}f(z)-b}),$$
which implies that
$$(\frac{f-a_{1}}{f-b_{1}})^{n}\equiv A(\frac{\Delta_{\eta}f-a_{2}}{\Delta_{\eta}f-b_{2}})^{m},$$
where $A$ is a nonzero constant. By Lemma 2.5, we obtain
\begin{eqnarray}
nT(r,f)=mT(r,\Delta_{\eta}f)+S(r,f).
\end{eqnarray}
It follows from (3.12) and (3.17) that $m=n$. Thus we get
\begin{eqnarray}
\frac{f-a}{f-b}\equiv C(\frac{\Delta_{\eta}f-a}{\Delta_{\eta}f-b}),
\end{eqnarray}
where $C$ is a nonzero constant. Since $(a_{2}-b_{2})f-(a_{1}-b_{1})\Delta_{\eta}f- a_{2}b_{1}+a_{1}b_{2}\not\equiv0$, hence $C\neq1$, then from above, we have
$$\frac{a_{2}-b_{2}}{\Delta_{\eta}f-a_{2}}\equiv \frac{(C-1)f-Cb_{1}+a_{1}}{f-a_{1}},$$
and
$$T(r,f)=T(r,\Delta_{\eta}f)+S(r,f).$$
Obviously, $\frac{Cb_{1}-a_{1}}{C-1}\neq a_{1}$ and $\frac{Cb_{1}-a_{1}}{C-1}\neq b_{1}$. It follows that $N(r,\frac{1}{f-\frac{Cb_{1}-a_{1}}{C-1}})=0$. Then by the Second Fundamental Theorem,
\begin{eqnarray*}
\begin{aligned}
2T(r,f)&\leq \overline{N}(r,f)+\overline{N}(r,\frac{1}{f-a_{1}})+\overline{N}(r,\frac{1}{f-b_{1}})+\overline{N}(r,\frac{1}{f-\frac{Cb_{1}-a_{1}}{C-1}})+S(r,f)\\
&\leq \overline{N}(r,\frac{1}{f-a_{1}})+\overline{N}(r,\frac{1}{f-b_{1}})+S(r,f),
\end{aligned}
\end{eqnarray*}
that is $2T(r,f)\leq \overline{N}(r,\frac{1}{f-a_{1}})+\overline{N}(r,\frac{1}{f-b_{1}})+S(r,f)$, which contradicts (3.1).

{\bf Subcase 1.2} \quad $H_{n,m}\not\equiv0$ for any positive integers $m$ and $n$. From the above discussion, we know that a zero of $f-a_{1}$ and $\Delta_{\eta}f-a_{2}$ (or a zero of $f-b_{1}$ and $\Delta_{\eta}f-b_{2}$ ) with multiplicities $m$ and $n$, must be the zero of $n\varphi-m\psi$.  So we have
\begin{align}
&\overline{N}_{(m,n)}(r,\frac{1}{f-a_{1}})+\overline{N}_{(m,n)}(r,\frac{1}{f-b_{1}})\leq \overline{N}(r,\frac{1}{n\varphi-m\psi})+S(r,f)\notag\\
&\leq T(r,n\varphi-m\psi)+S(r,f)=S(r,f).
\end{align}

Thus by (3.1), (3.12) and (3.19), we  get
\begin{align}
&T(r,f)\leq \overline{N}(r,\frac{1}{f-a_{1}})+\overline{N}(r,\frac{1}{f-b_{1}})+S(r,f)&\notag\\
&\leq N_{1}(r,\frac{1}{f-a_{1}})+\overline{N}_{[2]}(r,\frac{1}{f-a_{1}})+\overline{N}_{(3}(r,\frac{1}{f-a_{1}})\notag\\
&+N_{1}(r,\frac{1}{\Delta_{\eta}f-b_{2}})+\overline{N}_{[2]}(r,\frac{1}{\Delta_{\eta}f-b_{2}})+\overline{N}_{(3}(r,\frac{1}{\Delta_{\eta}f-b_{2}})+S(r,f)\notag\\
&\leq\frac{1}{15}(N(r,\frac{1}{f-a_{1}})+N(r,\frac{1}{f-b_{1}}))+\frac{1}{15}(N(r,\frac{1}{\Delta_{\eta}f-a_{2}})\notag\\
&+N(r,\frac{1}{\Delta_{\eta}f-b_{2}}))+\frac{2}{3}T(r,f)+S(r,f)\leq \frac{14}{15}T(r,f)+S(r,f),
\end{align}
it follows that $T(r,f)=S(r,f)$, a contradiction.

{\bf Case 2} \quad $\Delta_{\eta}\varphi\equiv0$. That is $\varphi$ is a periodic entire function with period $\eta$.

This completes Theorem 8.

\section{The proof of Theorem 9 }
Suppose on the contrary that such entire function $f$ of finite order exists. Since $f$ is an  entire function, and that $f$ and $\Delta_{\eta}f$ share $(0,0)$ and $(a,-a)$ IM, then  by Theorem 8, we know either $f_{\eta}=0$, which is impossible or $\Delta_{\eta}\varphi\equiv0$. That is 
\begin{align}
\frac{f'f_{\eta}}{f(f-a)}\equiv\frac{f'_{\eta}f_{2\eta}}{f_{\eta}(f_{\eta}+a)}.
\end{align}
Easy to see from (4.1) that $f_{\eta}$ and $\Delta_{\eta}f_{\eta}$ share $(0,0)$ and $(a,-a)$ IM. On the other hand, we can see from (3.1) that the zeros of $f_{\eta}$ almost consist of the zeros of $f-0$ and the zeros of $f-a$. Because $\varphi$ is a small function of $f$, so  almost all the zeros of $f_{\eta}$ and $f_{2\eta}$ are of multiplicity $1$.  That is to say, $f_{\eta}$ and $\Delta_{\eta}f_{\eta}$ share $0$ CM almost. Then there is a meromorphic function $h$ and an entire function $g$ such that
\begin{align}
\frac{\Delta_{\eta}f_{\eta}}{f_{\eta}}=he^{g}.
\end{align}
According to Lemma 2.1, we know that $T(r,f_{\eta})=T(r,\Delta_{\eta}f_{\eta})+S(r,f)$, which is (3.12) by Lemma 2.3. Similar to the proof of Case 1 of Theorem 8, we can obtain that $f_{\eta}\equiv0$, a contradiction.

This completes Theorem 9.

\section{The proof of Theorem 10 }
We prove by contradiction. Assume that $f(z)\not\equiv\Delta_{\eta}f(z)$. Then by Theorem 8, we only need to consider the Case that $\Delta_{\eta}\varphi\equiv0$. That is
\begin{align}
\varphi(z)\equiv \varphi(z+\eta).
\end{align}
We claim that $\varphi(z)\equiv C$, where $C$ is a constant. Otherwise, if $\varphi(z)\not\equiv C$, (3.22) can deduce that $\varphi(z)$ has zeros $z_{0}$, then easy to see that $z_{0}+\eta, z_{0}+2\eta,\ldots,$ all are zeros of $\varphi(z)$. Thus
\begin{align}
T(r,\varphi)\geq N(r,\frac{1}{\varphi})\geq (1+o(1))r,
\end{align}
and then
\begin{align}
\lambda(\varphi)\geq1.
\end{align}
From the fact that $f(z)$ and $\Delta_{\eta}f(z)$ share $a$ and $b$ IM, we obtain
\begin{align}
\lambda(f)=\lambda(\Delta_{\eta}f).
\end{align}
Easy to see from (3.7) that
\begin{eqnarray*}
\begin{aligned}
&2m(r,f(z))=m(r,\varphi(z)f^{2}(z) )+S(r,f)\\
&=m(r,f'(z)(f(z)-\Delta_{q,\eta}f(z))+(a+b)\varphi(z) f(z)+ab\varphi(z))+S(r,f)\\
&\leq m(r,f'(z)(f(z)-\Delta_{q,\eta}f(z))+(a+b)\varphi(z) f(z))+m(r,ab\varphi(z))+S(r,f)\\
&\leq m(r,f(z))+m(r,\frac{f'(z)(f(z)-\Delta_{\eta}f(z))}{f(z)})+m(r,(a+b)\varphi(z))+S(r,f)\\
&\leq m(r,f(z))+m(r,f'(z))+S(r,f),
\end{aligned}
\end{eqnarray*}
which deduces
\begin{align}
T(r,f(z))\leq T(r,f'(z))+S(r,f).
\end{align}
On the other hand
\begin{align}
T(r,f'(z))\leq T(r,f(z))+S(r,f).
\end{align}
Hence
\begin{align}
T(r,f(z))= T(r,f'(z))+S(r,f).
\end{align}
(3.28) implies
\begin{align}
\lambda(f)=\lambda(f')<1.
\end{align}
Then by Lemma 2.7, (5.3), (5.4) and (5.8), we have
\begin{align}
\lambda(\varphi)\leq\lambda(f)<1,
\end{align}
and then $\varphi(z)$ is of  order less than $1$, which contradicts with (5.3). And hence by Lemma 2.8, we can get that $\varphi(z)\equiv C$. If $C=0$, we can obtain that $f(z)\equiv\Delta_{\eta}f(z)$, a contradiction. So $C\neq0$. Choose the matrix $A$ to be 
\begin{equation} 
\left ( 
\begin{array}{ccc} 
0 & 4 & -1\\ 
-\frac{2}{3} & -8 & 2\\ 
\end{array}
\right).
\end{equation}
We can obtain that the rank of $A$ is $2$. Take $a_{3}=\frac{5}{4}$, $c_{1}=0$ and $c_{2}$ with $e^{c_{2}}=\frac{4}{9}a^{2}$. Apply these complex number to the equation (2.1), and by a tedious work, we can see that equation (2.1) is (3.2). Since the matrix $A$ we choose satisfies $a_{0}b_{1}-a_{1}b_{0}\neq0$, and $a_{1}b_{2}-a_{2}b_{1}\neq0$, then by Lemma 2.6, we can obtain that 
\begin{align}
f(z)=C_{1}e^{dz}+C_{2}e^{-dz}=e^{-dz}(C_{1}e^{2dz}+C_{2}).
\end{align}
Since $f(z)$ and $\Delta_{\eta}f(z)$ share $a$ and $-a$ IM, we can rewrite (3.2) as 
\begin{align}
\frac{1}{f}=\frac{f'(f-\Delta_{\eta}f)}{f(f-a)(f+a)\varphi}.
\end{align}
Then by (3.4) and  Lemma 2.1, we have
\begin{align}
m(r,\frac{1}{f})&=m(r,\frac{f'(f-\Delta_{\eta}f)}{f(f-a)(f+a)\varphi})\leq m(r,\frac{1}{\varphi})+m(r,\frac{f-\Delta_{\eta}f}{f})\notag\\
&+m(r,\frac{f-\Delta_{\eta}f}{f})+m(r,\frac{f'}{(f-a)(f+a)})=S(r,f).
\end{align}
Therefore, by the First Fundamental Theorem of Nevanlinna,  (5.11) and (5.13), we get
\begin{align}
T(r,f)=N(r,\frac{1}{f})+S(r,f)=N(r,\frac{e^{dz}}{C_{1}e^{2dz}+C_{2}})=2dT(r,e^{z})+S(r,f).
\end{align}

That is to say,  $\rho(f)=\rho(e^{z})=1$, it contradicts the assumption $\rho(f)<1$. From above discussion, we obtain $f(z)\equiv\Delta_{\eta}f(z)$.

This completes Theorem 10.

\

{\bf Acknowledgements} The author would like to thank to anonymous referees for their helpful comments.


\end{document}